\documentclass[reqno, 11pt, a4paper]{amsart}

\usepackage{fullpage}
\selectfont

\usepackage[utf8]{inputenc}
\usepackage[OT2,T1]{fontenc}
\usepackage[english]{babel}

\usepackage{amsmath}
\usepackage{amsfonts}
\usepackage{amssymb}
\usepackage{amsthm}

\usepackage{mathtools}
\usepackage{leftidx}
\usepackage{listings}
\usepackage{relsize}
\usepackage{bm}
\usepackage{tikz}
\usetikzlibrary{matrix}

\usepackage{breqn}
\newcounter{myequation}[equation]

\usepackage[backref=page,pagebackref=true,hyperindex=true,bookmarks=true]{hyperref}

\hypersetup{
  pdfauthor   = {Lercier, Ritzenthaler},
  pdftitle    = {Siegel modular forms of degree three and invariants of ternary quartics},
  pdfsubject  = {},
  pdfkeywords = {},
  colorlinks=true,
  breaklinks=true, urlcolor=blue, linkcolor=blue, citecolor=blue,
  bookmarksopen=true}

\theoremstyle{plain}
\newtheorem{theorem}{Theorem}[section]

\newtheorem{proposition}[theorem]{Proposition}

\theoremstyle{definition}
\newtheorem{definition}[theorem]{Definition}

\theoremstyle{remark}
\newtheorem{remark}[theorem]{Remark}

\numberwithin{equation}{section}

\DeclareMathOperator{\im}{Im}
\DeclareMathOperator{\Jac}{Jac}
\DeclareMathOperator{\SL}{SL}
\DeclareMathOperator{\Sp}{Sp}
\DeclareMathOperator{\tr}{tr}

\def\C{\mathbb{C}}
\def\H{\mathbb{H}}
\def\Q{\mathbb{Q}}
\def\R{\mathbb{R}}
\def\Z{\mathbb{Z}}

\def\alphab{\bm{\alpha}}
\def\betab{\bm{\beta}}
\def\gammab{\bm{\gamma}}
\def\deltab{\bm{\delta}}
\def\chib{\bm{\chi}}

\newcommand{\bfI}{\mathbf{I}}
\newcommand{\bfR}{\mathbf{R}}

\newcommand{\mtt}[1]{{\text{\relsize{-1.5}\tt\bf #1}}}

\newcommand{\car}[2]{%
  \left[{\substack{#1\\#2}}\right]}

\newcommand*\circled[1]{\tikz[baseline=(char.base)]{
    \node[shape=circle,draw,inner sep=1pt] (char) {#1};}}

\begin{document}

\title{Siegel modular forms of degree three\\ and invariants of ternary quartics}
\date{\today}

\begin{abstract}
  We determine the structure of the graded ring of Siegel modular forms of
  degree 3. It is generated by 19 modular forms, among which we identify a
  homogeneous system of parameters with 7 forms of weights $4$, $12$, $12$,
  $14$, $18$, $20$ and $30$.
  We also give a complete dictionary between the Dixmier-Ohno invariants of
  ternary quartics and the above generators.
\end{abstract}

\author[Lercier]{Reynald Lercier}
\address{%
  Reynald Lercier,
  DGA \& Univ Rennes, %
  CNRS, IRMAR - UMR 6625, F-35000
  Rennes, %
  France. %
}
\email{reynald.lercier@m4x.org}

\author[Ritzenthaler]{Christophe Ritzenthaler}
\address{%
  Christophe Ritzenthaler,
  Univ Rennes, %
  CNRS, IRMAR - UMR 6625, F-35000
  Rennes, %
  France. %
}
\email{christophe.ritzenthaler@univ-rennes1.fr}

\subjclass[2010]{14K20, 14K25, 14J15, 11F46, 14L24}
\keywords{Siegel modular forms, plane quartics, invariants, generators, explicit}

\maketitle

\section{Introduction and main results}
\label{sec:intr-main-results}

Let $g \geq 1$ be an integer and let $\bfR_g(\Gamma_g)$ denote the
$\C$-algebra of \emph{modular forms of degree $g$} for the symplectic group
$\Sp_{2g}(\Z)$ (see Section~\ref{sec:analyt-sieg-modul} for a precise
definition). It is a normal and integral domain of finite type over $\C$,
closely related to the moduli space of principally polarized abelian varieties
over $\C$. But even generators of these algebras are only known for small values of
$g$: $g=1$ is usually credited to Klein~\cite{klein,klein1897} and
Poincar\'e~\cite{poincare1905,poincare1911}, $g=2$ to Igusa
\cite{igusa-thetag2} and $g=3$ to Tsuyumine \cite{Tsuyumine86}. In the latter,
Tsuyumine gives $34$ generators and asks if they form a minimal set of
generators.  We answer  in the negative and  prove in the present paper that there exists a subset of $19$
of them which still generates the algebra and which is minimal
(Theorem~\ref{th:mfbasis}). As a by-product we also exhibit a (possibly
incomplete) set of $55$ relations and use them to obtain a homogeneous system
of parameters for this algebra (Theorem~\ref{th:hsop}).

Unlike Tsuyumine, we extensively use computer algebra software since we base our
strategy on evaluation/interpolation which leads to computing ranks and invert
large dimensional matrices. Still, a naive application of this strategy would have forced us to
work with complex numbers, which would have been bad for efficiency but also
to certify our computations. Hence, in order to perform exact arithmetic
computations, we make a detour through the beautiful geometry of smooth plane
quartics and Weber's formula \cite{weber} which allows us to express values
(of quotients) of the theta constants and ultimately modular forms as rational
numbers (up to a fourth root of unity). The strategy could be interesting for
future investigations for $g=4$ as those theta constants can be computed in a
similar way \cite{turku}.

We then move on to a second task in the continuation of the famous Klein's
formula, see \cite[Eq. 118, p. 462]{klein} and \cite{LRZ, matone,
  ichi-klein}. This formula relates a certain modular form of weight $18$,
namely $\chib_{18}$, to the square of the discriminant of plane quartics. A
complete dictionary between modular forms and invariants was only known for
$g=1$ and $g=2$. For $g=3$, these formulas can come in two flavors:
restricting to the the image of the hyperelliptic locus in the Jacobian locus,
one gets expressions of the modular forms in terms of Shioda invariants for
binary octics, see~\cite{Tsuyumine86} and~\cite{lorenzo19}; considering the
generic case, one gets expressions in terms of Dixmier-Ohno invariants for
ternary quartics, see Proposition~\ref{prop:phiimg}. Extra care was taken in
making these formulas as normalized as possible using the background of
\cite{LRZ} and also to eliminate parasite coefficients
coming from relations between the invariants  as much as possible.
As a striking example, the modular form $\chib_{28}$ is equal to
$-2^{171} \cdot 3^3\:I_{27}^3\,I_3$ (the exponent of $2$ is large because the
normalization chosen by Dixmier for $I_{27}$ is not optimal at $2$). We
finally give formulas in the opposite direction and express all Dixmier-Ohno
invariants as quotients of modular forms by powers of $I_{27}$, see
Proposition~\ref{prop:phiinv}. We hope that such formulas may eventually lead
to a set of generators for the ring of invariants of ternary quartics with
good arithmetic properties.  Indeed, theta constants have intrinsically good
``reduction properties modulo primes'' (in the sense that they often have a
primitive Fourier expansion) and may help guessing such a set of generators.

The full list of expressions for the $19$ Siegel modular forms either in terms
of the theta constants or in terms of curve invariants, the expressions of
Dixmier-Ohno invariants in terms of Siegel modular forms and the $55$
relations in the algebra, are available at~\cite{LR2019}.

\subsection*{Acknowledgments} We warmly thank the anonymous referees for
carefully reading this work and for suggestions. This work is partially
supported by the French National Research Agency under the
\textsc{anr}-{\footnotesize 18}-\textsc{ce}{\footnotesize 40}-{\footnotesize
  0026}-{\footnotesize 01} \textsc{cl{\footnotesize ap}-cl{\footnotesize ap}}
project.

\section{Review of Tsuyumine's construction of  Siegel modular forms}
\label{sec:analyt-sieg-modul}

We recall here the definition of the $34$ generators for the $\C$-algebra of
modular forms of degree $3$ built by Tsuyumine. Surprisingly, they all are
 polynomials in theta constants with rational coefficients: one knows that
when $g \geq 5$, there exists modular forms which are not in the algebra
generated by theta constants \cite{manni}, while the answer for $g=4$ is
pending \cite{oura}. We take special care of the multiplicative constant
involved in each expression.

\subsection{Theta functions and theta constants}
Let $g \geq 1$ be an integer and $\H_g=\{\tau \in \mathbf{M}_{g}(\C), \; {^t}{\tau} = \tau, \ \im \tau > 0\}$.
\begin{definition}
  The \emph{theta function} with characteristics
  $\car{\varepsilon_{1}}{\varepsilon_{2}} \in \mathbf{M}_{2,g}(\Z)$ is given, for $z \in \C^g$ and
  $\tau \in \H_{g}$, by
  \begin{displaymath}
    \theta\,{\car{\varepsilon_{1}}{\varepsilon_{2}}}(z, \tau) =
    \sum_{n\, \in\, \Z^g}
    \exp(\,i \pi\; (n+\varepsilon_{1}/2) \:  \tau \: \ltrans{\,(n+\varepsilon_{1}/2)\,)}\ %
    \exp(\,2 i \pi\; (n + \varepsilon_{1}/2) \: \ltrans{\,(z + \varepsilon_{2}/2)\,)}.
  \end{displaymath}
The \emph{theta constant} (with characteristic $\car{\varepsilon_{1}}{\varepsilon_{2}}$)  is the function of $\tau$ defined as  $\theta\,{\car{\varepsilon_{1}}{\varepsilon_{2}}}(\tau)=    \theta\,{\car{\varepsilon_{1}}{\varepsilon_{2}}}(0, \tau)$.
\end{definition}

\begin{proposition}
  Let $z \in \C^g$,
  $\tau \in \H_{g}$, $\car{\varepsilon_{1}}{\varepsilon_{2}} \in
  \mathbf{M}_{2,g}(\Z)$, then
  \begin{equation}\label{eq:1}
    \theta\,{\car{\varepsilon_1}{\varepsilon_2}}(-z, \tau) =%
    \theta\,{\car{-\varepsilon_1}{-\varepsilon_2}}(z, \tau)\,,
  \end{equation}
  and
  \begin{equation}\label{eq:2}
    \forall \car{\delta_{1}}{\delta_{2}} \in
    \mathbf{M}_{2,g}(2\,\Z),\ \ \ \theta\,{\car{\varepsilon_1+\delta_1}{\varepsilon_2+\delta_2}}(z, \tau)
    = %
    \exp(i \pi\; \varepsilon_1\ltrans{\,\delta_2}/2) \:\: %
    \theta\,{\car{\varepsilon_1}{\varepsilon_2}}(z, \tau)\,.%
  \end{equation}
\end{proposition}

Combining these two equations shows that
$z \mapsto \theta\,{\car{\varepsilon_1}{\varepsilon_2}}(z,\tau)$ is even if
$\varepsilon_1\: \ltrans{\,\varepsilon_2} \equiv 0 \pmod{2}$, and odd
otherwise. The characteristics $\car{\varepsilon_{1}}{\varepsilon_{2}}$ are
then said to be even and odd, respectively.

The modular group
$\Gamma_g:=\Sp_{2g}(\Z)$ acts on $\H_{g}$ by
\begin{equation}\label{eq:9}
  \tau \rightarrow M.\tau := (\mtt{A}\tau+\mtt{B}) \, (\mtt{C}\tau+\mtt{D})^{-1}\quad \text{for} \
  M = \begin{psmallmatrix} \mtt{A} & \mtt{B} \\ \mtt{C} & \mtt{D} \end{psmallmatrix} \,,
\end{equation}
and on characteristics by
  \begin{displaymath}
    \car{\varepsilon_{1}}{\varepsilon_{2}} \rightarrow
    M.\car{\varepsilon_{1}}{\varepsilon_{2}} = %
    (\varepsilon_{1}\,^\frown\, \varepsilon_{2})\, M +
    (\ltrans{\,\mtt{A}}\,\mtt{C})_\Delta\,^\frown\,
    (\ltrans{\,\mtt{B}}\,\mtt{D})_\Delta.
  \end{displaymath}
    Here, ``$\,^\frown$'' denotes the concatenation of two row vectors, and
  ``$(.)_\Delta$'' denotes the row vector equal to the diagonal of the square
  matrix given in argument.
These result in the following action of $\Gamma_g$ on theta constants.

\begin{proposition}[{Transformation formula \cite[Chap.~5,
    Th.~2]{Igusa72}\cite[p.442]{manni1989}~\cite[Prop.~3.1.24]{Cosset11}}]
  \label{prop:theta-transform}
  Let $\tau \in \H_{g}$,
  $\car{\varepsilon_{1}}{\varepsilon_{2}} \in \mathbf{M}_{2,g}(\R)$ and
  $M\in\Gamma_g$, then
  \begin{equation}\label{eq:3}
    \theta\,{\car{\varepsilon_{1}}{\varepsilon_{2}}}(M.\tau) = %
    \zeta_M\,\sqrt{\det(\mtt{C}\tau+\mtt{D}
      )}\,\exp(\,-i\pi\;\sigma /
    4)\:\:\theta\,{\car{\delta_{1}}{\delta_{2}}}(\tau)
  \end{equation}
  with
  $\car{\delta_{1}}{\delta_{2}} = M.\car{\varepsilon_{1}}{\varepsilon_{2}}$\,,\:
  $\zeta_M$ an eighth root of unity depending only on $M$ and
  \begin{displaymath}
    \sigma= %
    \varepsilon_{1} \; \mtt{A} \, \ltrans{\,\mtt{B}}\;\ltrans{\,\varepsilon_{1}}%
    +2\; \varepsilon_{1}\;\mtt{B} \, \ltrans{\,\mtt{C}}\;\ltrans{\,\varepsilon_{2}}%
    + \varepsilon_{2} \; \mtt{C} \, \ltrans{\,\mtt{D}}\;\ltrans{\,\varepsilon_{2}}%
    + (\,2\; \varepsilon_{1}\;\mtt{A} + 2\; \varepsilon_{2}\;\mtt{C} +
    (\ltrans{\,\mtt{A}}\,\mtt{C})_\Delta) \;
    \ltrans{\,(\ltrans{\,\mtt{B}}\,\mtt{D})_\Delta}\,.
  \end{displaymath}
\end{proposition}

In the following, we only make use of theta constants with
characteristics with coefficients in $\{0,1\}$. Using Eq.~\eqref{eq:2} in
combinaison with Eq.~\eqref{eq:3} allows to have a transformation formula purely between characteristics of this
form.

To lighten notations, we number the theta constants as in \cite{KLLRSS17}. We write
\begin{math}
  \theta_n :=
  \theta\,{\car{\delta_0\:\delta_1\:\ldots\:\delta_{g-1}}{\varepsilon_0\:\varepsilon_1\:\ldots\:\varepsilon_{g-1}}}
\end{math}
where $0 \leq n < 2^{2g-1}$ is the integer whose binary expansion is
``$\delta_0\cdots\delta_{g-1}\varepsilon_0\ldots\varepsilon_{g-1}$''.
In genus 3, there are 36 even theta constants (the odd ones are all $0$). We give
in Table~\ref{tab:tsuthetas} the correspondence between their numbering
 in~\cite[pp.789--790]{Tsuyumine86} and our binary numbering.

\begin{table}[htbp]
  \centering
  \renewcommand{\arraystretch}{1.2}
  \setlength{\arraycolsep}{2pt}
  $\begin{footnotesize}\begin{array}{l|cccccccccccccccccc|}\cline{2-19}%
    \text{Tsuyumine} & \circled{1} & \circled{2} & \circled{3} & \circled{4} & \circled{5} &%
    \circled{6} & \circled{7} & \circled{8} & \circled{9} & \circled{10} &%
    \circled{11} & \circled{12} & \circled{13} & \circled{14} & \circled{15} &%
    \circled{16} & \circled{17} &\circled{18} \\\hline
    \text{Binary} & \theta_{31} &\theta_{27} &\theta_{56} &\theta_{48} &\theta_{49} &\theta_{59} &\theta_{24} &\theta_{16} &\theta_{17} &\theta_{28} &\theta_{20} &\theta_{21} &\theta_{55} &\theta_{54} &\theta_{62} &\theta_{47} &\theta_{12} &\theta_{4}\\\hline
  \multicolumn{19}{c}{}\\[-10pt]\cline{2-19}%
    \text{Tsuyumine} & \circled{19} &\circled{20}&%
    \circled{21} &\circled{22} & \circled{23} & \circled{24} & \circled{25} &%
    \circled{26} & \circled{27} & \circled{28} & \circled{29} & \circled{30} &%
    \circled{31} & \circled{32} & \circled{33} &\circled{34} &\circled{35}&%
    \circled{36}\\\hline
    \text{Binary} & \theta_{5} &\theta_{8} &\theta_{0} &\theta_{1} &\theta_{35} &\theta_{34} &\theta_{42} &\theta_{40} &\theta_{32} &\theta_{33} &\theta_{3} &\theta_{2} &\theta_{10} &\theta_{7} &\theta_{6} &\theta_{14} &\theta_{45} &\theta_{61}\\\hline
  \end{array}\end{footnotesize}$
  \caption{Tsuyumine's numbering of even theta constants}
  \label{tab:tsuthetas}
\end{table}

\subsection{Siegel modular forms}
\label{sec:modular-groups}

 Let $\Gamma_g(\ell)$ denote the principal congruence subgroup of level $\ell$,
\textit{i.e.}
\begin{math}
  \{ M \in \Gamma_g\:|\: M \equiv \mathbf{1}_{2g} \bmod \ell\},
\end{math}
and let $\Gamma_g(\ell,2\ell)$ denote the congruence subgroup
\begin{math}
  \{ M \in \Gamma_g(\ell)\:|\: %
  (\ltrans{\,\mtt{A}}\,\mtt{C})_\Delta  \equiv
  (\ltrans{\,\mtt{B}}\,\mtt{D})_\Delta \equiv 0 \bmod 2\,\ell\}.
\end{math}

 For a congruence subgroup $\Gamma \subset \Gamma_g$, let
 $\mathbf{R}_{g,h}(\Gamma)$ be the $\C$-vector space of \emph{analytic Siegel
   modular forms of weight $h$ and degree $g$ for $\Gamma$}, consisting of
 complex holomorphic functions $f$ on $\H_{g}$ satisfying
\begin{displaymath}
f(M.\tau) = \det (\mtt{C} \tau+\mtt{D})^{h} \cdot f(\tau)
\end{displaymath}
for all
 $M \in \Gamma$.
For $g=1$, one also requires that $f$ is holomorphic at ``infinity'' but we
will not look at this case here.  We also denote the \emph{$\C$-algebra of
  Siegel modular forms of degree $g$ for $\Gamma$} by
$\mathbf{R}_g(\Gamma):= \bigoplus \mathbf{R}_{g,h}(\Gamma)$\,.
The modular group acts on $\mathbf{R}_{g,h}(\Gamma_g)$ by
\begin{displaymath}
  f \rightarrow M.f:=\det(\mtt{C} \tau+\mtt{D})^{-h} \cdot f(M.\tau)\,.
\end{displaymath}
In particular, $f \in  \mathbf{R}_{g,h}(\Gamma)$ if and only if $M.f=f$ for all $M \in \Gamma$.

We now restrict to $g=3$. A strategy to build modular forms for $\Gamma_3$ is first to construct a form
$F \in \mathbf{R}_3(\Gamma_3(2))$, and then average over the finite quotient
$\Gamma_3/\Gamma_3(2)$ to get a modular form $f \in \mathbf{R}_3(\Gamma_3)$,
namely
\begin{equation}\label{eq:4}
  f = \sum_{M\in \Gamma_3/\Gamma_3(2)} M.F\;.
\end{equation}
All forms $F$ which will be considered are polynomials in the theta constants, and are of even
weight. Hence, given an $F$, a careful application of the transformation
formula (Proposition.~\ref{prop:theta-transform}) gives all summands, where we do not
care about the choice of the square root as it is raised to an even power.

Tsuyumine gives some of the building blocks $F$s in terms of maximal syzygetic sets
of even characteristics~\cite[Sec.~21]{Tsuyumine86}. Multiplying the theta
constants in a given set is an element of $\mathbf{R}_3(\Gamma(2))$. The
quotient $\Gamma_3/\Gamma_3(2)$ acts transitively on these $135$ sets numbered
from $((1))$ to $((135))$ by Tsuyumine. Among them, $33$ are actually used to
define a set of generators for $\mathbf{R}_3(\Gamma_3)$. We give their
expressions in Table~\ref{tab:syzy}.

\begin{table}[htbp]
  \centering
  \begin{footnotesize}
    \begin{tabular}{cc}
      \setlength{\arraycolsep}{2pt}
      $\begin{array}{l|l}
         \multicolumn{1}{c|}{\#} & \multicolumn{1}{c}{\theta\text{-monomial}} \\\hline
         ((1)) & \ \,\theta_{3}\,\theta_{28}\,\theta_{31}\,\theta_{33}\,\theta_{34}\,\theta_{61}\,\theta_{62}\,\theta_{0} \\
         ((2)) & \text{-}\,\theta_{1}\,\theta_{2}\,\theta_{28}\,\theta_{31}\,\theta_{32}\,\theta_{35}\,\theta_{61}\,\theta_{62} \\
         ((3)) & \text{-}\,\theta_{3}\,\theta_{8}\,\theta_{20}\,\theta_{31}\,\theta_{33}\,\theta_{42}\,\theta_{54}\,\theta_{61} \\
         ((4)) & \ \,\theta_{1}\,\theta_{10}\,\theta_{20}\,\theta_{31}\,\theta_{35}\,\theta_{40}\,\theta_{54}\,\theta_{61} \\
         ((5)) & \ \,\theta_{2}\,\theta_{8}\,\theta_{21}\,\theta_{31}\,\theta_{32}\,\theta_{42}\,\theta_{55}\,\theta_{61} \\
         ((18)) & \text{-}\,\theta_{1}\,\theta_{4}\,\theta_{17}\,\theta_{20}\,\theta_{40}\,\theta_{45}\,\theta_{56}\,\theta_{61} \\
         ((31)) & \ \,\theta_{3}\,\theta_{4}\,\theta_{7}\,\theta_{24}\,\theta_{27}\,\theta_{28}\,\theta_{31}\,\theta_{0} \\
         ((32)) & \text{-}\,\theta_{1}\,\theta_{2}\,\theta_{5}\,\theta_{6}\,\theta_{24}\,\theta_{27}\,\theta_{28}\,\theta_{31} \\
         ((34)) & \ \,\theta_{1}\,\theta_{5}\,\theta_{10}\,\theta_{14}\,\theta_{16}\,\theta_{20}\,\theta_{27}\,\theta_{31} \\
         ((36)) & \text{-}\,\theta_{4}\,\theta_{10}\,\theta_{14}\,\theta_{17}\,\theta_{21}\,\theta_{27}\,\theta_{31}\,\theta_{0} \\
         ((37)) & \text{-}\,\theta_{4}\,\theta_{7}\,\theta_{28}\,\theta_{31}\,\theta_{32}\,\theta_{35}\,\theta_{56}\,\theta_{59} \\\hline
       \end{array}$
    &
      $\begin{array}{l|l}
         \multicolumn{1}{c|}{\#} & \multicolumn{1}{c}{\theta\text{-monomial}} \\\hline
         ((38)) & \ \,\theta_{5}\,\theta_{6}\,\theta_{28}\,\theta_{31}\,\theta_{33}\,\theta_{34}\,\theta_{56}\,\theta_{59} \\
         ((39)) & \text{-}\,\theta_{7}\,\theta_{8}\,\theta_{16}\,\theta_{31}\,\theta_{32}\,\theta_{47}\,\theta_{55}\,\theta_{56} \\
         ((43)) & \ \,\theta_{7}\,\theta_{12}\,\theta_{20}\,\theta_{31}\,\theta_{35}\,\theta_{40}\,\theta_{48}\,\theta_{59} \\
         ((45)) & \ \,\theta_{7}\,\theta_{24}\,\theta_{31}\,\theta_{40}\,\theta_{47}\,\theta_{48}\,\theta_{55}\,\theta_{0} \\
         ((47)) & \text{-}\,\theta_{14}\,\theta_{17}\,\theta_{31}\,\theta_{33}\,\theta_{47}\,\theta_{48}\,\theta_{62}\,\theta_{0} \\
         ((51)) & \text{-}\,\theta_{1}\,\theta_{6}\,\theta_{24}\,\theta_{31}\,\theta_{40}\,\theta_{47}\,\theta_{49}\,\theta_{54} \\
         ((54)) & \text{-}\,\theta_{3}\,\theta_{12}\,\theta_{16}\,\theta_{31}\,\theta_{34}\,\theta_{45}\,\theta_{49}\,\theta_{62} \\
         ((55)) & \ \,\theta_{3}\,\theta_{24}\,\theta_{27}\,\theta_{32}\,\theta_{35}\,\theta_{56}\,\theta_{59}\,\theta_{0} \\
         ((73)) & \ \,\theta_{8}\,\theta_{16}\,\theta_{24}\,\theta_{32}\,\theta_{40}\,\theta_{48}\,\theta_{56}\,\theta_{0} \\
         ((85)) & \ \,\theta_{1}\,\theta_{16}\,\theta_{17}\,\theta_{32}\,\theta_{33}\,\theta_{48}\,\theta_{49}\,\theta_{0} \\
         ((89)) & \text{-}\,\theta_{2}\,\theta_{3}\,\theta_{4}\,\theta_{5}\,\theta_{48}\,\theta_{49}\,\theta_{54}\,\theta_{55} \\\hline
       \end{array}$\\&\\
    \multicolumn{2}{c}{
      $\begin{array}{l|l}
         \multicolumn{1}{c|}{\#} & \multicolumn{1}{c}{\theta\text{-monomial}} \\\hline
         ((90)) & \ \,\theta_{1}\,\theta_{6}\,\theta_{7}\,\theta_{48}\,\theta_{49}\,\theta_{54}\,\theta_{55}\,\theta_{0} \\
         ((99)) & \text{-}\,\theta_{5}\,\theta_{8}\,\theta_{17}\,\theta_{28}\,\theta_{34}\,\theta_{47}\,\theta_{54}\,\theta_{59} \\
         ((103)) & \ \,\theta_{4}\,\theta_{8}\,\theta_{12}\,\theta_{16}\,\theta_{20}\,\theta_{24}\,\theta_{28}\,\theta_{0} \\
         ((111)) & \ \,\theta_{1}\,\theta_{4}\,\theta_{5}\,\theta_{16}\,\theta_{17}\,\theta_{20}\,\theta_{21}\,\theta_{0} \\
         ((115)) & \text{-}\,\theta_{8}\,\theta_{20}\,\theta_{28}\,\theta_{34}\,\theta_{42}\,\theta_{54}\,\theta_{62}\,\theta_{0} \\
         ((118)) & \text{-}\,\theta_{3}\,\theta_{10}\,\theta_{21}\,\theta_{28}\,\theta_{33}\,\theta_{40}\,\theta_{55}\,\theta_{62} \\
         ((119)) & \text{-}\,\theta_{1}\,\theta_{20}\,\theta_{21}\,\theta_{34}\,\theta_{35}\,\theta_{54}\,\theta_{55}\,\theta_{0} \\
         ((131)) & \ \,\theta_{1}\,\theta_{2}\,\theta_{3}\,\theta_{4}\,\theta_{5}\,\theta_{6}\,\theta_{7}\,\theta_{0} \\
         ((132)) & \text{-}\,\theta_{4}\,\theta_{5}\,\theta_{6}\,\theta_{7}\,\theta_{32}\,\theta_{33}\,\theta_{34}\,\theta_{35} \\
         ((133)) & \ \,\theta_{2}\,\theta_{8}\,\theta_{10}\,\theta_{32}\,\theta_{34}\,\theta_{40}\,\theta_{42}\,\theta_{0} \\
         ((135)) & \ \,\theta_{1}\,\theta_{2}\,\theta_{3}\,\theta_{32}\,\theta_{33}\,\theta_{34}\,\theta_{35}\,\theta_{0} \\\hline
       \end{array}$
      }
    \end{tabular}
  \end{footnotesize}\medskip
  \caption{Tsuyumine's maximal syzygetic sequences}
  \label{tab:syzy}
\end{table}

Then Tsuyumine considers $34$  $F$s written as combinations of
\begin{itemize}
\item $\chi_{18}=\prod_{\theta_i\text{ even}}\theta_i$,
\item a rational function of the 36 non-zero $\theta_i^{\,4}$\,,
\item the monomials $((i))$ defined in Table~\ref{tab:syzy}\,, and
\item the squares of the $\gcd$ between two such $((i))$.
\end{itemize}

\begin{table}[htbp]
  \begin{footnotesize}
    \renewcommand{\arraystretch}{1.1}
    \setlength{\arraycolsep}{2pt}
    \begin{math}\hspace*{-0.6cm}\begin{array}{l|c|c|r}
    \multicolumn{1}{c|}{\text{Name}} &%
    \multicolumn{1}{c|}{\text{[Tsu86] Coeff.}} &%
    \multicolumn{1}{c|}{F\in\mathbf{R}_3(\Gamma(2))} &%
    \multicolumn{1}{c}{\#\text{sum.}}\\\hline\hline%
    \chi_{18} & \ 1/(2^9 \cdot 3^4\cdot 5\cdot  7)                  & \prod_{\theta_i\text{ even}}\theta_i & 1\\
    \chi_{28} & \ 1/(2^{10}\cdot 3^2\cdot 5\cdot 7)        & \chi_{18}^2\,/\,((131))^2 & 135\\\hline\hline
    \alpha_4 & \ 1/(2^{13}\cdot 3\cdot 7)            & \gcd(\:((131)), ((132))\:)^2 & 945 \\
    \alpha_6 & \ 1/(2^6\cdot 3\cdot 7)             & \theta_{0}^4  \cdot  ((131)) & 1080 \\
    \alpha_{10} & \text{-}1/(2^4\cdot 3^2\cdot 5\cdot 11)       & (\theta_{16}\,\theta_{20}\,\theta_{32}\,\theta_{34}\,\theta_{48}\,\theta_{54})^2 \cdot ((131)) & 30240\\
    \alpha_{12} & \ 1/(2^8\cdot 3^5\cdot 5)           & (\theta_{2}\,\theta_{21}\,\theta_{24}\,\theta_{49}\,\theta_{62}\,\theta_{0})^4 & 336 \\
    \alpha_{12}' & \ 3/2^8                   &  ((85))^2 \cdot ((119))^2\,/\,(\,\theta_{1}\,\theta_{0}\,)^4 & 945 \\
    \alpha_{16} & \text{-}3^2/2^9              & ((85))^2 \cdot ((119)) \cdot ((131)) & 3780\\
    \alpha_{18} & \text{-}3^2/2^5              & \theta_{0}^4\,((85))^2 \cdot ((119)) \cdot ((131)) &7560\\
    \alpha_{20} & \ 3/(2^9\cdot 5)               & (((85))^2 \cdot ((119))^2 \cdot ((131))^2\,/\,(\theta_{1}\,\theta_{0})^4 & 63\\
    \alpha_{24} & \ 3^2/2^3                 & \theta_{0}^4\,((85))^2 \cdot ((119))^2 \cdot ((131))^2\,/\,\theta_{1}^4 & 1260 \\
    \alpha_{30} & \ 3^4/(2^8\cdot 5)             & (((85))^3 \cdot ((119))^3 \cdot ((131))^3\,/\,(\theta_{1}\,\theta_{0}^2)^4 & 1260 \\\hline\hline
    \beta_{14}  & \ 1/(2^5\cdot 3\cdot 7)             & \theta_{31}^8\, \chi_{18}/(\,((5)) \cdot ((54))\,) & 4320 \\
    \beta_{16}  & \ 1/(2^6\cdot 3)               & ((31)) \cdot ((43)) \cdot ((47)) \cdot ((51)) & 7560 \\
    \beta_{22}  & \text{-}1/(2^4\cdot 3)              & (\theta_{27}\,\theta_{31}\,\theta_{54}\,\theta_{55}\,\theta_{59}\,\theta_{62})^4\,\chi_{18}\,/\,(\,((2)) \cdot ((54))\,) & 30240 \\
    \beta_{22}' & \ 2^4                     & \chi_{18}\,((119))^2 \cdot ((133))^2\,/\,(\,\theta_{34}^4\theta_{0}^4\,((18)) \cdot ((34))\,) & 90720 \\
    \beta_{26} & \text{-}1/2^2                  & ((32)) \cdot ((36)) \cdot ((37)) \cdot ((45)) \cdot ((90)) \cdot ((111)) \cdot ((135))\,/\,\theta_{0}^4 & 362880 \\
    \beta_{28} & \text{-}1/2^2                  & ((32)) \cdot ((36)) \cdot ((37)) \cdot ((45)) \cdot ((90)) \cdot ((111)) \cdot ((135)) & 362880 \\

    \beta_{32} & \ 1/2^2                   & \chi_{18}\,((85))^2 \cdot  ((89))^2 \cdot  ((90)) \cdot  ((111)) \cdot ((135))\,/\,(\,\theta_{48}^4\,\theta_{49}^4\,\theta_{0}^4 ((4)) \cdot  ((99))\,) & 362880 \\

    \beta_{34} & \ 1/(2^3\cdot 3)               & \theta_{31}^8\,\chi_{18}\, ((90))^2 \cdot  ((111))^2 \cdot  ((135))^2\,/\,(\,\theta_{0}^4\, \theta_{1}^4\, ((3)) \cdot  ((31))\,) & 120960 \\\hline\hline
   \gamma_{20} & \ 1/(2^7\cdot 3)               & \theta_{31}^4\,\chi_{18}\,((135))\,/\,((1)) & 7560 \\
   \gamma_{24} & \ 1/2^7                   & \theta_{31}^8\,\chi_{18}^2\,/\,(\,((4)) \cdot ((5)) \cdot ((47)) \cdot ((54)) \,) & 11340 \\
   \gamma_{26} & \ 1/2^6                   & (\theta_{31}\,\theta_{28})^4\,\chi_{18}\,((38)) \cdot ((135))\,/\,((1)) & 22680 \\
   \gamma_{32} & \ 1/(2^3\cdot 3)               & (\theta_{16}\,\theta_{20}\,\theta_{31}\,\theta_{49}\,\theta_{54}\,\theta_{56}\,\theta_{59})^4\,\chi_{18}\,((135))\,/\,((1)) & 120960 \\
  c_{32}' & \text{-}1/2^3                  & \theta_{33}^4\,\chi_{18}\,((90))^2 \cdot ((111))^2 \cdot ((135))\,/\,(\,(\theta_{1}\,\theta_{0})^4\,((1))\,) & 30240 \\
   \gamma_{36} & \text{-}1/2^4                  & (\theta_{28}\,\theta_{31})^4\,\chi_{18}\,((38)) \cdot ((90)) \cdot ((111)) \cdot ((135))^2\,/\,(\,\theta_{1}^4\,((1))\,) & 181440 \\
   \gamma_{38} & \ 1/2^4                   & \theta_{31}^{16}\,\chi_{18}^2\,((31)) \cdot ((39)) \cdot ((43))\,/\,(\,\theta_{7}^4\,((4)) \cdot ((5)) \cdot ((47)) \cdot ((54))\,) & 90720 \\
       c_{38}' & \ 1/2^2                   & \theta_{31}^4\,\chi_{18}\,((38))^2 \cdot ((90)) \cdot ((111)) \cdot ((135))^2\,/\,(\,\theta_{1}^4\,((1))\,) & 362880 \\
   \gamma_{42} & \ 1/2^3                   & \chi_{18}\,(((38)) \cdot ((85))^2 \cdot ((90)) \cdot ((111)) \cdot ((119))^2 \cdot ((135))\,/\,(\,(\theta_{1}\,\theta_{0})^4\,((1))\,) & 181440 \\
   \gamma_{44} & \ 1/2^4                   & \chi_{18}^2\,\theta_{31}^8\,((45))^2 \cdot ((55))^2 \cdot ((103))^2\,/\,(\,(\theta_{24}\,\theta_{0})^4\,((4)) \cdot ((5)) \cdot ((47)) \cdot ((54))\,) &90720 \\\hline\hline
   \delta_{30} & \ 2^7/3                   & (\theta_{28}\,\theta_{31})^4\,\chi_{18}\, ((47)) \cdot ((115)) \cdot ((118))\,/\,((1)) & 90720 \\
   \delta_{36} & \ 1/2^3                   & (\theta_{28}\,\theta_{31}\,\theta_{0})^4\,\chi_{18}\,((31)) \cdot ((38)) \cdot ((118)) \cdot ((135))\,/\,((1)) & 181440 \\
   \delta_{46} & \text{-}1/2                    & (\theta_{28}\,\theta_{31})^4\,\chi_{18}\,((31)) \cdot ((38)) \cdot ((90)) \cdot ((111)) \cdot ((118)) \cdot ((135))^2\,/\,((1)) & 725760 \\
   c_{48} & \ 1/2                     & \theta_{28}^4\,\chi_{18}\,((31))^2 \cdot ((38)) \cdot ((90)) \cdot ((111)) \cdot ((118)) \cdot ((135))^2\,/\,((1)) & 725760 \\\hline
  \end{array}\end{math}
  \end{footnotesize}\medskip
  \caption{Tsuyumine's generators (the index is their weight), Tsuyumine's
    normalization constant, the form $F$ and
    the number of summands of the polynomial in the theta constants}
  \label{tab:tsubasis}
\end{table}

Using the map from modular forms to invariants of binary octics introduced by
Igusa~\cite{igusag3}, he proves the following result.
\begin{theorem}[Tsuyumine~\mbox{\cite[Sec.~20]{Tsuyumine86}}\footnote{See~\cite[p.~44]{tsushima89}
  for the $(1-T^{12})$ misprint in the denominator of Equation~(\ref{eq:5})
  in~\cite{Tsuyumine86}.}]\label{th:tsuyumine}
The graded algebra $\mathbf{R}_3(\Gamma_3)$ is generated by the 34 modular
forms defined in Table~\ref{tab:tsubasis}. Its Hilbert–Poincaré series is
generated by the rational function
\begin{equation}\label{eq:5}
  \frac{(\, 1+{T}^{2} \,) \: N(T)}{(1-T^4)\,(1-T^{12})^2\,(1-T^{14})\,(1-T^{18})\,(1-T^{20})\,(1-T^{30})}\,,
\end{equation}
where
  \begin{dgroup*}[style={\small},spread={-2pt}]
    \begin{dmath*}
      N(T) = 1 - {T}^{2} + {T}^{4} + {T}^{10} + 3\,{T}^{16} - {T}^{18} +
      3\,{T}^{20} + 2\,{T}^{22} + 2\,{T}^{24} + 3\,{T}^{26} + 4\,{T}^{28} +
      2\,{T}^{30} + 7\,{T}^{32} + 3\,{T}^{34} + 7\,{T}^{36} + 5\,{T}^{38} +
      9\,{T}^{40} + 6\,{T}^{42} + 10\,{T}^{44} + 8\,{T}^{46} + 10\,{T}^{48} +
      9\,{T}^{50} + 12\,{T}^{52} + 7\,{T}^{54} + 14\,{T}^{56} + 7\,{T}^{58} +
      12\,{T}^{60} + 9\,{T}^{62} + 10\,{T}^{64} + 8\,{T}^{66} + 10\,{T}^{68} +
      6\,{T}^{70} + 9\,{T}^{72} + 5\,{T}^{74} + 7\,{T}^{76} + 3\,{T}^{78} +
      7\,{T}^{80} + 2\,{T}^{82} + 4\,{T}^{84} + 3\,{T}^{86} + 2\,{T}^{88} +
      2\,{T}^{90} + 3\,{T}^{92} - {T}^{94} + 3\,{T}^{96} + {T}^{102} +
      {T}^{108} - {T}^{110} + {T}^{112}\,.
    \end{dmath*}
  \end{dgroup*}
\end{theorem}

The modular forms $f$ defined in Table~\ref{tab:tsubasis} are all polynomials
in the theta constants whose  primitive part has all its coefficients
equal to $\pm 1$ and whose content is
\begin{displaymath}
c(f) = \frac{\# \Gamma_3/\Gamma_3(2)}{\#\{\text{summands of } f\}}=
\frac{2^9 \cdot 3^4 \cdot 5 \cdot 7}{\#\{\text{summands of } f\}} \in \Z\,.
\end{displaymath}
In order to get simpler expressions when restricting to the hyperelliptic
locus or to the decomposable one, Tsuyumine multiplies each $f$ by an
additional normalization constant ($2^{\mathrm{nd}}$ column of
Table~\ref{tab:tsubasis}). For instance, as defined by Tsuyumine,
  \begin{displaymath}
    \chi_{28} := 2^{-10}\cdot 3^{-2}\cdot 5^{-1}\cdot 7^{-1} \sum_{M\in
      \Gamma_g/\Gamma_g(2)} M.(\chi_{18}^2\,/\,((131))^2)\,,
  \end{displaymath}
  and therefore the 135 summands are each a (monic) monomial in the theta constants
  times
  $\pm (2^{-10}\cdot 3^{-2}\cdot 5^{-1}\cdot 7^{-1}) \cdot c(\chi_{28})= \pm
  1/30$ (the sign depends on the monomial).

  Having in mind possible applications of our results to fields of positive
  characteristic, we replace the multiplication by Tsuyumine's constant by a
  multiplication by $1/c(f)$. In this way, $f$ is a sum of (monic) monomials
  in the theta constants with coefficients $\pm 1$. To avoid confusion with
  Tsuyumine's notation, our modular forms will be denoted with bold
  font. Typically, $\chib_{28} := 30\,\chi_{28}$,
  $\alphab_4 := 112\,\alpha_4$, $\alphab_6 := \alpha_6$,
  $\alphab_{10} := 165\,\alpha_{10}$, etc.

  Still driven by the link with the hyperelliptic locus, Tsuyumine adds to
  $c_{32}'$ (resp. $c_{38}'$ and $c_{48}$) some polynomials in modular forms
  of smaller weights and denote the result $\gamma_{32}$ (resp. $\gamma_{38}$
  and $\delta_{48}$). Theorem~\ref{th:tsuyumine} as stated in
  \cite{Tsuyumine86} considers modular forms $\gamma_{32}$, $\gamma_{38}$ and
  $\delta_{48}$, instead of $c_{32}'$, $c_{38}'$ and $c_{48}$. The two
  theorems are obviously equivalent. Here, we choose instead to define
  $\gammab_{32}:=c_{32}'/6$, $\gammab_{38}:=c_{38}'$ and
  $\deltab_{48}:=c_{48}$\,.

  \begin{remark} \label{rem:eva} Some of the modular forms in
    Table~\ref{tab:tsubasis} have a large number of summands. While it would be
    cumbersome to store them, evaluating them is relatively quick as it
    basically consists in permuting theta constants up to some eighth roots of
    unity according to Eq.~(\ref{eq:3}).  Following Tsuyumine, the sum is
    computed in two steps. Let $\Theta$ be the  subgroup of $\Gamma_3$ conjugate to
    $\Gamma_3(1,2)$ that stabilizes $\theta_{61}$ ($\Gamma_3(1,2)$ stabilizes
    $\theta_{0}$).  Tsuyumine gives explicit coset representatives for
    $\Gamma_3/\Theta$ (36 elements) and $\Theta/\Gamma_3(2)$ ($8!$ elements)
    and splits the sum in Eq.~(\ref{eq:4}) as
\begin{displaymath}
  f = \sum_{M'\in \Gamma_3/\Theta} M'\,.\,\sum_{M''\in \Theta/\Gamma_3(2)}  M''.F
\end{displaymath}
\noindent
We use this approach in order to perform the computation of the
summands\footnote{There are two small typos
  in~\cite[pp. 842--846]{Tsuyumine86}, the $(3,6)$-th coefficients of
  ``$M_{\scalebox{0.5}{\circled{1}}}$'' must be -1 instead of 1, and the
  $(2,2)$-th coefficients of ``$M_{\scalebox{0.5}{\circled{27}}}$'' must be 1
  instead of $0$. This modification makes $M_{\scalebox{0.5}{\circled{1}}}$
  and $M_{\scalebox{0.5}{\circled{27}}}$ symplectic.}.  In order to do that, we
also need the eighth roots of unity $\zeta_{M'}$ and $\zeta_{M''}$ from Proposition.~\ref{prop:theta-transform}.  One
approach is to precompute them using a fixed chosen matrix in $\H_3$.
A better solution is, with the notation of Eq.~(\ref{eq:9}), to make use of
the relation $\zeta_M^4=(-1)^{\tr(\mtt{B}\:
  \ltrans{\mtt{C}})}$~\cite[Chap. 5]{Igusa72}. Since the modular forms have
even weight, the degree of $F$ in the theta constants is a multiple of 4, as
well as the powers of $\zeta_{M'}$ and $\zeta_{M''}$.
\end{remark}

\section{A minimal set of generators for modular forms of degree $3$}
\label{sec:newbase}

\subsection{Fundamental set of modular forms}
\label{sec:fund-modul-forms}

Since we know the dimensions of each $\mathbf{R}_{3,h}(\Gamma_3)$ from the
generating functions of Theorem~\ref{th:tsuyumine}, it is a matter of linear
algebra to check that a given subset of Tsuyumine's generators is enough for
generating the full algebra.
It would involve choosing a monomial ordering on the ring of theta constants and
computing a Gr\"obner basis of the homogeneous ideal defined by the generating
subset given as formal expressions in terms of them (see~\cite[Section 1.4.1]{DK2015}).
However, it is difficult to perform these computations since there
exist numerous algebraic relations between the theta constants.  Therefore we
favor an interpolation/evaluation strategy as follows.

Suppose that we want to prove that a given form $f$ of weight $h$, given as a polynomial in the theta constants, can be
obtained from a given set $\{f_1,\ldots,f_m\}$. This set produces
$F_1,\ldots,F_n$, homogeneous polynomials in the $f_i$ of weight $h$.  If
$n<d=\dim \mathbf{R}_{3,h}(\Gamma_3)$, then all forms of weight $h$ cannot be
obtained. Assume that $n \geq d$. Then, if we can find
$(\tau_i)_{i=1,\dots,d} \in \H_g^d$ such that the matrix
$(F_i(\tau_j))_{1 \leq i,j \leq d}$ is of rank $d$, we know that $f$ can be
written in terms of the $f_i$, and even find such a relation. Equivalently, we
will actually find a polynomial relation between $f/\theta_{0}^{2h}$ and the
$f_i/\theta_{0}^{2 w_i}$ where $w_i$ denotes the weight of $f_i$.

By Remark~\ref{rem:eva}, the evaluation of a form
$f(\tau)/\theta_{0}^{2h}(\tau)$ boils down to the computation of quotients
$(\theta_i/\theta_{0})(\tau)$.  A naive approach would be to use
an arbitrary matrix  $\tau \in\H_3$. But then the theta constants would in general be
transcendental complex numbers which would make the computations much more costly and the final result hard to certify.
We therefore prefer to consider
a complex torus $\Jac C$ attached to a smooth plane quartic $C$ given by an
Aronhold system. Indeed (see for instance \cite{weber, agmri, NR17}), $7$
general lines in $\P^2$ form an Aronhold system of $7$ bitangents for a unique
plane quartic $C$. Then, one can easily recover the equations of the $21$
other bitangents and an expression of the quotients
$ {(\theta_i/\theta_{0})^4(\tau)}$ in terms of the coefficients of the linear
forms defining the bitangents (see for instance \cite[Theorems~2 and
3]{NR17}). Note that we do not explicitly know the Riemann matrix $\tau$ here,
since it depends not only on $C$ but also on the choice of a symplectic basis
for $H_1(C,\Z)$. But when each of the bitangents in the Aronhold system is
defined over $\Q$, all computations can be performed over $\Q$ and
$ {(\theta_i/\theta_{0})^4(\tau)}$ is a rational number.

To remove the fourth root of unity ambiguity that remains, we start by
computing independently an approximation over $\C$ of an explicit Riemann
matrix $\tau'$ for the curve $C$. We need to do it only at very low precision
(a typical choice is 20 decimal digits) and this can be done efficiently
either in \textsc{maple} (package \textsc{algcurves} by Deconinck et
al.~\cite{MR1837895}) or in \textsc{magma} (package \textsc{riemann surfaces}
by Neurohr~\cite{neurohr18}). Then, we can calculate an approximation of the
theta constants at $\tau'$.

To conclude, note that \cite[Theorem~3.1]{NR17} shows that
the set $\left\{ {\theta_j^{\,8}} / {\theta_i^{\,8}} \right\}$ running through
every even theta constants $\theta_i,\theta_j$ depends only on $C$ and not on
the Riemann matrix. Indeed, the dependence on this matrix relies only on the
quadratic form $q_{_0}$ (in the notation of \textit{loc. cit.}) whose
contribution disappears in the eighth power. Therefore, there exist
an integer $i_0$ and a permutation $\sigma$ such that
$$\frac{\theta_{\sigma(i)}(\tau')^{\,8}}{\theta_{i_0}(\tau')^{\,8}}=\frac{\theta_{i}(\tau)^{\,8}}{\theta_{0}(\tau)^{\,8}}.$$
We simply enumerate all the possible candidates for $i_0$ until we find a
suitable $\sigma$ that gives $i_0$ and $\sigma$.
Then, since we know ${\theta_{\sigma(i)}(\tau')}/{\theta_{i_0}(\tau')}$ with
small precision and its eighth power exactly, it is possible to obtain the
exact value of ${\theta_{i}(\tau)}/{\theta_{0}(\tau)}$.\medskip

Using this method extensively  leads to a set of generators for
$\mathbf{R}_{3}(\Gamma_3)$. Moreover it is easy to prove, by the same
algorithms, that this set is \emph{fundamental}, \textit{i.e.} one cannot remove any element and still generate the algebra
$\mathbf{R}_{3}(\Gamma_3)$.
\begin{theorem}\label{th:mfbasis}
  The $19$ Siegel modular forms $\alphab_{4}$, $\alphab_{6}$, $\alphab_{10}$,
  $\alphab_{12}$, $\alphab_{12}'$, $\betab_{14}$, $\alphab_{16}$,
  $\betab_{16}$, $\chib_{18}$, $\alphab_{18}$, $\alphab_{20}$, $\gammab_{20}$,
  $\betab_{22}$, $\betab_{22}'$, $\alphab_{24}$, $\gammab_{24}$,
  $\gammab_{26}$, $\chib_{28}$ and $\alphab_{30}$ define a fundamental set of generators for
  $\mathbf{R}_{3}(\Gamma_3)$.
\end{theorem}
\begin{remark}
Note that \cite{runge2} proved that $\mathbf{R}_{3}(\Gamma_3(2))$ has a fundamental set of  generators of $30$ elements.
\end{remark}

A word on the complexity. The proof mainly consists in checking for all the
even weight $h$ between 4 and 48 that there exists an evaluation matrix of
rank $\dim \mathbf{R}_{3,h}(\Gamma_3)$ for this set of $19$ modular forms. It
is a matter of few hours for the largest weight to perform this calculation in
\textsc{magma}. Most of the time is spent on the evaluation of the $19$ forms
$f_i$ at a matrix $\tau_j$, which takes about 1~minute on a laptop.\smallskip

Additionally, we find the expressions of the remaining 15 modular forms given in
Table~\ref{tab:tsubasis}. The first
ones are\medskip
\begin{dgroup*}[style={\footnotesize},spread={-2pt}]
  \begin{dmath*}
    2^5\cdot 3^4\cdot 5\cdot 7^2 \cdot 11\:\: \betab_{26} =%
7\,\alphab_{6}\,\alphab_{10}^2 - 3080\,\alphab_{6}^2\,\betab_{14} - 145530\,\alphab_{12}\,\betab_{14} + 194040\,\alphab_{12}'\,\betab_{14} - 11760\,\alphab_{10}\,\alphab_{16} - 7040\,\alphab_{4}\,\alphab_{6}\,\betab_{16} + 16660\,\alphab_{10}\,\betab_{16} -
    20824320\,\alphab_{4}^2\,\chib_{18} - 4435200\,\alphab_{6}\,\alphab_{20} + 2822512\,\alphab_{6}\,\gammab_{20} - 55440\,\alphab_{4}\,\betab_{22} +
    36960\,\alphab_{4}\,\betab_{22}' - 105557760\,\gammab_{26}\,,
  \end{dmath*}
  \begin{dsuspend}\end{dsuspend}
  \begin{dmath*}
    2^8\cdot3^4\cdot5^2\cdot7^4\:\: \betab_{28} =  %
    -105\,\alphab_{4}^2\,\alphab_{10}^2 - 42000\,\alphab_{4}^2\,\alphab_{6}\,\betab_{14} + 66885\,\alphab_{4}\,\alphab_{10}\,\betab_{14} + 129654\,\betab_{14}^2 - 96000\,\alphab_{4}^3\,\betab_{16} + 77792400\,\alphab_{12}\,\betab_{16} +
    207446400\,\alphab_{12}'\,\betab_{16} + 5399533440\,\alphab_{4}\,\alphab_{6}\,\chib_{18} - 9996323400\,\alphab_{10}\,\chib_{18} - 4321800\,\alphab_{10}\,\alphab_{18} + 320544000\,\alphab_{4}^2\,\alphab_{20} +
    82576256\,\alphab_{4}^2\,\gammab_{20} - 12965400\,\alphab_{6}\,\betab_{22} - 17287200\,\alphab_{6}\,\betab_{22}' - 666792000\,\alphab_{4}\,\alphab_{24} - 700378560\,\alphab_{4}\,\gammab_{24} -
    442172001600\,\chib_{28}\,,
  \end{dmath*}
  \begin{dsuspend}\end{dsuspend}
  \begin{dmath*}
    2^3\cdot3^4\cdot5\cdot7^4\:\: \deltab_{30} = %
    -37044\,\betab_{14}\,\betab_{16} + 23040\,\alphab_{4}^3\,\chib_{18} +
    987840\,\alphab_{6}^2\,\chib_{18} + 47508930\,\alphab_{12}\,\chib_{18} +
    133358400\,\alphab_{12}'\,\chib_{18} - 1568\,\alphab_{4}\,\alphab_{6}\,\gammab_{20}
    + 46305\,\alphab_{10}\,\gammab_{20} - 246960\,\alphab_{6}\,\gammab_{24} +
    282240\,\alphab_{4}\,\gammab_{26}\,,
  \end{dmath*}
  \begin{dsuspend}\end{dsuspend}
  \begin{dmath*}
    2\cdot3\cdot5\cdot7\:\: \gammab_{32}' = %
    \chib_{18}\,(\alphab_{4}\,\alphab_{10}  - 252\,\betab_{14})\,.
  \end{dmath*}
\end{dgroup*}
The last ones, for instance $\gammab_{44}$, $\deltab_{46}$ and $\deltab_{48}$,
tend to be heavily altered with the relations that exist between these 19
modular forms, and have huge coefficients (thousands of digits).

\subsection{Module of relations between the generators}
\label{sec:relations}

We now quickly deal with the relations defining the algebra
$\mathbf{R}_{3}(\Gamma_3)$. With the same techniques, involving modular forms
up to weight $70$ (see Remark~\ref{rem:modforminDO} for speeding up the
computations), we find a (possibly incomplete) list of $55$ relations for our
generators of $\mathbf{R}_{3}(\Gamma_3)$ given by weighted polynomials of
degree $32$ to $58$ (\textit{cf.} Table~\ref{tab:rels}).

\begin{table}[htbp]
  \centering
  \begin{tabular}{l|c|c|c|c|c|c|c|c|c|c|c|c|c|c|}
    {Weight} & 32 & 34 & 36 & 38 & 40 & 42 & 44 & 46 & 48 & 50 & 52 & 54 &56 & 58 \\\hline
    {Number} & 1 & 1 & 2 & 4 & 5 & 5 & 7 & 6 & 8 & 6 & 5 &  2 & 2 & 1 \\\hline
  \end{tabular}
  \caption{number of relations of a given weight in   $\mathbf{R}_{3}(\Gamma_3)$}
  \label{tab:rels}
\end{table}

The relations of weight $32$ and $34$ are relatively small,
\begin{dgroup*}[style={\footnotesize},spread={-2pt}]
  \begin{dmath*}
    0 =%
    -25226544365568\,\betab_{16}^2 + 50854572195840\,\betab_{16}\,\alphab_{16}
    - 25092716544000\,\alphab_{16}^2 +
    13916002383360\,\alphab_{18}\,\betab_{14} -
    18410871153185280\,\chib_{18}\,\betab_{14} +
    1109304189987840\,\gammab_{20}\,\alphab_{12p} -
    1951854879744000\,\alphab_{20}\,\alphab_{12p} -
    413549252645760\,\gammab_{20}\,\alphab_{12} +
    1463891159808000\,\alphab_{20}\,\alphab_{12} +
    474409172160\,\betab_{22p}\,\alphab_{10} +
    355806879120\,\betab_{22}\,\alphab_{10} +
    8471592360\,\alphab_{12p}\,\alphab_{10}^2 -
    3882813165\,\alphab_{12}\,\alphab_{10}^2 +
    14993672601600\,\gammab_{26}\,\alphab_{6} -
    1800579432960\,\betab_{14}\,\alphab_{12p}\,\alphab_{6} +
    559752621120\,\betab_{14}\,\alphab_{12}\,\alphab_{6} +
    14755739264\,\betab_{16}\,\alphab_{10}\,\alphab_{6} -
    25299240960\,\alphab_{16}\,\alphab_{10}\,\alphab_{6} -
    4775514472960\,\gammab_{20}\,\alphab_{6}^2 +
    10174277836800\,\alphab_{20}\,\alphab_{6}^2 -
    43285228\,\alphab_{10}^2\,\alphab_{6}^2 +
    7065470720\,\betab_{14}\,\alphab_{6}^3 +
    779296133468160\,\chib_{28}\,\alphab_{4} -
    530133424128\,\betab_{14}^2\,\alphab_{4} -
    2857212610560\,\betab_{16}\,\alphab_{12p}\,\alphab_{4} +
    1510363895040\,\betab_{16}\,\alphab_{12}\,\alphab_{4} -
    5020202880\,\alphab_{18}\,\alphab_{10}\,\alphab_{4} +
    59052646477440\,\chib_{18}\,\alphab_{10}\,\alphab_{4} -
    104866460160\,\betab_{22p}\,\alphab_{6}\,\alphab_{4} -
    38488222080\,\betab_{22}\,\alphab_{6}\,\alphab_{4} +
    16149647360\,\betab_{16}\,\alphab_{6}^2\,\alphab_{4} +
    642585968640\,\gammab_{24}\,\alphab_{4}^2 +
    516363724800\,\alphab_{24}\,\alphab_{4}^2 +
    1529966592\,\betab_{14}\,\alphab_{10}\,\alphab_{4}^2 -
    5609877504000\,\chib_{18}\,\alphab_{6}\,\alphab_{4}^2 -
    130817347584\,\gammab_{20}\,\alphab_{4}^3 -
    154557849600\,\alphab_{20}\,\alphab_{4}^3 -
    1036728\,\alphab_{10}^2\,\alphab_{4}^3 +
    97574400\,\betab_{14}\,\alphab_{6}\,\alphab_{4}^3 +
    223027200\,\betab_{16}\,\alphab_{4}^4\,,
\end{dmath*}
\begin{dsuspend}\end{dsuspend}
\begin{dmath*}
  0 =%
  -113265734400\,\alphab_{18}\,\betab_{16}
  -107036119008000\,\chib_{18}\,\betab_{16}
  +130691232000\,\alphab_{18}\,\alphab_{16}
  +123503214240000\,\chib_{18}\,\alphab_{16}
  +711613758240\,\gammab_{20}\,\betab_{14}
  +242595599400\,\betab_{22}\,\alphab_{12}'
  -121297799700\,\betab_{22}'\,\alphab_{12}
  -107820266400\,\gammab_{24}\,\alphab_{10}
  +670881657600\,\chib_{28}\,\alphab_{6}
  -399334320\,\betab_{14}^2\,\alphab_{6}
  +2662228800\,\betab_{16}\,\alphab_{12}'\,\alphab_{6}
  +3993343200\,\betab_{16}\,\alphab_{12}\,\alphab_{6}
  +80673600\,\alphab_{18}\,\alphab_{10}\,\alphab_{6}
  +699198091200\,\chib_{18}\,\alphab_{10}\,\alphab_{6}
  -221852400\,\betab_{22}'\,\alphab_{6}^2
  -665557200\,\betab_{22}\,\alphab_{6}^2
  +657308736\,\betab_{16}\,\betab_{14}\,\alphab_{4}
  -4978713600\,\alphab_{16}\,\betab_{14}\,\alphab_{4}
  +37811907302400\,\chib_{18}\,\alphab_{12}'\,\alphab_{4}
  +16298463535200\,\chib_{18}\,\alphab_{12}\,\alphab_{4}
  +5427686880\,\gammab_{20}\,\alphab_{10}\,\alphab_{4}
  +21254365440\,\gammab_{24}\,\alphab_{6}\,\alphab_{4}
  -2545060\,\betab_{14}\,\alphab_{10}\,\alphab_{6}\,\alphab_{4}
  -194295615360\,\chib_{18}\,\alphab_{6}^2\,\alphab_{4}
  -2123573760\,\gammab_{26}\,\alphab_{4}^2
  +27165600\,\betab_{14}\,\alphab_{12}'\,\alphab_{4}^2
  -50935500\,\betab_{14}\,\alphab_{12}\,\alphab_{4}^2
  -7299040\,\betab_{16}\,\alphab_{10}\,\alphab_{4}^2
  +9466800\,\alphab_{16}\,\alphab_{10}\,\alphab_{4}^2
  -1339238208\,\gammab_{20}\,\alphab_{6}\,\alphab_{4}^2
  +5145\,\alphab_{10}^2\,\alphab_{6}\,\alphab_{4}^2
  +5174400\,\betab_{22}'\,\alphab_{4}^3
  -19404000\,\betab_{22}\,\alphab_{4}^3
  -4011279360\,\chib_{18}\,\alphab_{4}^4\,.
\end{dmath*}
\end{dgroup*}

Runge~\cite[Cor.6.3]{runge1} shows that $\bfR_3(\Gamma_3)$ is a Cohen-Macaulay
algebra. There exists a strong link between a minimal free resolution of a
Cohen-Macaulay algebra and its Hilbert series. Let us rewrite
Equation~(\ref{eq:5}) as a rational fraction with denominator
$\prod_{d_i}(1 - T^{d_i})$ where the degrees $d_i$ run through the weights of
the fundamental set of generators.  We obtain a numerator with 140 non-zero
coefficients, the first and last ones of which are
\begin{footnotesize}
  \begin{multline*}
    1 - T^{32} - T^{34} - 2\,T^{36} - 4\,T^{38} - 5\,T^{40} - 5\,T^{42} -
    7\,T^{44} - 6\,T^{46} - 8\,T^{48} - 5\,T^{50} - 4\,T^{52} \\ + 4\,T^{56} +
    9\,T^{58} + 15\,T^{60} + 22\,T^{62} + 27\,T^{64} + 32\,T^{66} +
    36\,T^{68} + 39\,T^{70} + 36\,T^{72} + 34\,T^{74} + 26\,T^{76} + \ldots\\
    \ldots-5\,T^{296} - 8\,T^{298} - 6\,T^{300} - 7\,T^{302} - 5\,T^{304} -
    5\,T^{306} - 4\,T^{308} - 2\,T^{310} - T^{312} - T^{314} + T^{346}.
  \end{multline*}
\end{footnotesize}
The coefficients of the numerator give information on the weights and numbers
of relations.  They are consistent with Table~\ref{tab:rels} up to weight
$48$. The drop from $6$ (relations) to a coefficient $5$ in weight $50$
indicates that there is a first syzygy (\textit{i.e.} a relation between the
relations) of weight~50.

\subsection{A homogeneous system of parameters}
\label{sec:homog-syst-param}

Having these relations, one can also try to work out a homogeneous system of
parameters (\textsc{hsop}) for $\mathbf{R}_{3}(\Gamma_3)$. Recall that this is
a set of elements $(f_i)_{1\leq i\leq m}$ of the algebra, which are algebraically
independent, and such that $\mathbf{R}_{3}(\Gamma_3)$ is a
$\C[f_1,\ldots,f_m]$-module of finite type. Equation~(\ref{eq:5}) suggests
that a \textsc{hsop} of weight $4$, $12$, $12$, $14$, $18$, $20$ and $30$ may
exist. An easy Gr\"obner basis computation made in \textsc{magma} with the
lexicographic order
$\alphab_{6} < \alphab_{10}< \ldots < \gammab_{26}<\chib_{28}$ shows that when
we set to zero $\alphab_4$, $\alphab_{12}$, $\alphab_{12}'$, $\betab_{14}$,
$\chib_{18}$, $\alphab_{20}$ and $\alphab_{30}$ in the $55$ relations of
Table~\ref{tab:rels}, the remaining $12$ Siegel modular forms of the generating
set of Theorem~\ref{th:mfbasis} must be zero as well. As it is well known that
the dimension of
$\textrm{Proj}(\mathbf{R}_{3}(\Gamma_3))$ is $6$,
this yields the following theorem.

\begin{theorem}\label{th:hsop}
  A homogeneous system of parameters for $\mathbf{R}_{3}(\Gamma_3)$ is given
  by the $7$ forms $\alphab_4$, $\alphab_{12}$, $\alphab_{12}'$,
  $\betab_{14}$, $\chib_{18}$, $\alphab_{20}$ and $\alphab_{30}$.
\end{theorem}

\section{A dictionary between modular forms and invariants of quartics}
\label{sec:dictionary}

In~\cite{dixmier}, Dixmier gives a homogeneous system of parameters for the
graded $\C$-algebra $\bfI_3$ of invariants of ternary quartic forms under the
action of $\SL_3(\C)$. They are denoted $I_3$, $I_6$, $I_9$, $I_{12}$,
$I_{15}$, $I_{18}$ and $I_{27}$.  This list is completed by Ohno with six
invariants, $J_9$, $J_{12}$, $J_{15}$, $J_{18}$, $I_{21}$ and $J_{21}$, into a list of
13 generators for $\bfI_3$, the so-called \emph{Dixmier-Ohno
  invariants}~\cite{ohno,elsenhans}. Note that $2^{40} \cdot I_{27} = D_{27}$
where $D_{27}$ denotes the normalized discriminant of plane quartics in the
sense of \cite[p.426]{gelfand} or \cite[Prop.11]{demazure}.

Using the morphism $\rho_3$ defined in \cite{igusag3}, Tsuyumine in
\cite[pp. 847--864]{Tsuyumine86} relates each of the Siegel modular forms
given in Table~\ref{tab:tsubasis} with an invariant for the graded ring of
binary octics under the action of $\SL_2(\C)$. He uses this key argument to
prove Theorem~\ref{th:tsuyumine}. More generally, there is a way to
canonically associate an invariant to a modular form. After briefly recalling
the way to do so when $g=3$, we establish a complete dictionary between
$\bfR_3(\Gamma_3)$ and $\bfI_3$.

\subsection{Modular forms in terms of invariants}
\label{sec:modular-forms-terms}

Let us recall from \cite[2.2]{LRZ} how to associate an element of $\bfI_3$ to
$f \in \bfR_{3,h}(\Gamma_3)$. This morphism only depends on
the choice of a universal basis of regular differentials $\omega$ which can be
fixed ``canonically'' for smooth plane quartics (in the sense that it is a
basis of regular differentials over $\Z$).  Let $Q \in \C[x_1,x_2,x_3]$ be a ternary quartic form such
that $C: Q=0$ is a smooth genus $3$ curve.
Let $\Omega = \left[\begin{smallmatrix}\Omega_{1} \\
    \Omega_{2}\end{smallmatrix}\right]$ be the $6 \times 3$ period matrix of
$C$ defined by integrating $\omega_C$
with respect to an arbitrary symplectic basis of
$H_1(C,\Z)$. We have $\tau = \Omega_2^{-1} \Omega_1 \in \H_3$.  The function
\begin{equation}\label{eq:6}
  Q \mapsto \Phi_3(f)(Q) = \left(\frac{(2i\pi)^{3}}{\det \Omega_2}\right)^h
  \cdot f(\tau)
\end{equation}
is a homogeneous element of $\bfI_3$ of degree $3h$ (identifying the polynomial with its polynomial function).
\begin{remark}
  A similar construction can be worked out with invariants of binary octics
  (see \cite{ionica}). Up to a normalization constant, this is actually the
  same morphism as defined by~\cite{igusag3}.
\end{remark}
Chai's expansion principle \cite{chai} shows that if the Fourier expansion of
$f$ has coefficients in a ring $R \subset \C$, then $\Phi_3(f)$ is defined
over $R$ as well. When $f$ is given by a polynomial in the theta constants
with coefficients in $\Z$, we can take $R=\Z$.  A particular case is given by
the modular form $\chib_{18}$ which is the product of the $36$ theta
constants.  In \cite{LRZ} (see also \cite{ichi-klein}) one shows the following
precise form of Klein's formula~\cite[Eq. 118, p. 462]{klein},
\begin{equation}\label{eq:7}
\Phi_3(\chib_{18})=-2^{28} \cdot D_{27}^2 = -2^{28} \cdot
(\,2^{40}\,I_{27}\,)^2\,.
\end{equation}

\begin{remark}\label{rem:mod11}
  The map \eqref{eq:6} is obtained by pulling back \emph{geometric modular
    forms} to invariants as described in \cite{LRZ}. Within this background,
  it is for instance possible to speak about the reduction modulo a prime of
  modular forms and to consider the algebra that they generate. In small
  characteristics, one still encounters similar accidents as in the case of
  invariants. We will not study this question further here, but for instance,
  our $19$ generators have a surprising congruence modulo $11$,
  \begin{displaymath}
      \betab_{16} + 9\,\alphab_{16} + 3\,\alphab_{10}\,\alphab_{6} = 0 \bmod 11\,.
  \end{displaymath}
\end{remark}\medskip

We have seen in Section~\ref{sec:newbase} that we have an
evaluation/interpolation strategy to handle quotient of modular forms by a
power of $\theta_{0}$. This strategy can also be used to find the relations
with invariants. But now, we also need to take care of the transcendental
factor $\mu:= {(2i\pi)^{3}}/{\det \Omega_2}$. This is done in the following
way.
\begin{enumerate}
\item Assume that a relation $\Phi_3(f_0)=I_0$ is known  for a
  modular form $f_0$ of weight $h_0$. This is the case for $\chib_{18}$
  (\textit{cf.} Eq.~(\ref{eq:7})) and we will start with this one, but switch
  to a relation of lower weight  (\textit{i.e.} $4$ with $\alphab_4$ or even $2$ with
  $\chib_{18}/\alphab_4^4$\,) after a first round of the following steps (this
  simplifies the last step).
\item Let now $f$ be one of the generators from Theorem~\ref{th:mfbasis} of
  weight $h$ and compute a basis $j_1,\ldots,j_d$ of invariants of degree
  $3h$. We aim at finding $a_1,\ldots,a_d \in \Q$ such that
  $\Phi_3(f)=\sum a_i j_i$. This is done by evaluation/interpolation at
  Riemann models until one gets a system of $d$ linearly independent
  equations. More precisely, for a given $Q=0$ and an associated $\tau \in \H_3$\,:
\begin{enumerate}
\item Compute the values of $(j_1,\ldots,j_d)$ at $Q$;
\item Using the same procedure as in Section~\ref{sec:newbase}, compute
  $(f/\theta_{0}^{2h})(\tau)$ and $(f_0/\theta_{0}^{2h_0})(\tau)$;
\item Let $p=\textrm{lcm}(h_0,h)$. Since
  \begin{displaymath}
    \frac{(f/\theta_{0}^{2h})^{p/h}}{(f_0/\theta_{0}^{2h_0})^{p/h_0}} = \frac{(\mu^h f)^{p/h}}{(\mu^{h_0} f_0)^{p/h_0}} = \frac{\Phi_3(f)^{p/h}}{\Phi_3(f_0)^{p/h_0}}\,,
  \end{displaymath}
  we get the value of $\Phi_3(f)^{p/h}$. An approximate computation at low
  precision can then give the exact value.
\end{enumerate}
\end{enumerate}
The above strategy provides explicit expressions for $\Phi_3(f)$ where $f$ is
any modular form in the fundamental set defined in Theorem~\ref{th:mfbasis}.

\begin{proposition}\label{prop:phiimg}
  Let $f$ be a modular form of weight $h$ from
  Theorem~\ref{th:mfbasis}. There exists an explicit polynomial $P_f$ of
  degree $3h$ in the Dixmier-Ohno invariants such that
  \begin{displaymath}
    \Phi_3(f) = P_f(\, I_3, I_6, \ldots, I_{27} \,)\,.
  \end{displaymath}
\end{proposition}
\noindent
The first ones\footnote{We make available the list of these 19 polynomials at~\cite[file
``\texttt{SiegelMFfromDO.txt}'']{LR2019}.} are
\begin{dgroup*}[style={\footnotesize},spread={-2pt}]
  \begin{dmath*}
    \Phi_3(\alphab_4) =\ %
    2^{20} \cdot 3^3\cdot 7\:\:(%
    486\, I_{12} - 155520\, I_{6}^2 - 423\, J_{9}\, I_{3} + 117\, I_{9}\,
    I_{3} + 14418\, I_{6}\, I_{3}^2 + 8\, I_{3}^4%
    )\,,
  \end{dmath*}
  \begin{dmath*}
    5 \cdot 7\:\:
    \Phi_3(\alphab_6) = %
    \text{-}\,2^{28} \cdot 3^4\:\:(%
    40415760\,J_{18} -1224720\,I_{18} -2664900\,J_{9}^2 -8323560\,J_{9}\,I_{9}
    +2506140\,I_{9}^2 -76982400\,J_{12}\,I_{6} -1143538560\,I_{12}\,I_{6}
    +135992908800\,I_{6}^3 -40041540\,J_{15}\,I_{3} +2143260\,I_{15}\,I_{3}
    +247160160\,J_{9}\,I_{6}\,I_{3} +289325520\,I_{9}\,I_{6}\,I_{3}
    +400950\,J_{12}\,I_{3}^2 -6206220\,I_{12}\,I_{3}^2
    -7357573440\,I_{6}^2\,I_{3}^2 +1527453\,J_{9}\,I_{3}^3
    -266481\,I_{9}\,I_{3}^3 -36764280\,I_{6}\,I_{3}^4 -62720\,I_{3}^6%
    )\,,
  \end{dmath*}
  \begin{dmath*}
    \Phi_3(\alphab_{12}) =\ %
    2^{75}\cdot 3 \:\:(%
    495\,I_{27}\,J_{9} - 261\,I_{27}\,I_{9} - 14580\,I_{27}\,I_{6}\,I_{3} +
    32\,I_{27}\,I_{3}^3%
    )\,,
  \end{dmath*}
  \begin{dmath*}
    \Phi_3(\betab_{14}) =\ %
    2^{81}\cdot 3^4 \:\:(%
    -540\,I_{27}\,J_{15} - 4860\,I_{27}\,I_{15} + 285120\,I_{27}\,J_{9}\,I_{6}
    - 45360\,I_{27}\,I_{9}\,I_{6} + 810\,I_{27}\,J_{12}\,I_{3} +
    12204\,I_{27}\,I_{12}\,I_{3} - 18057600\,I_{27}\,I_{6}^2\,I_{3} -
    8541\,I_{27}\,J_{9}\,I_{3}^2 + 2961\,I_{27}\,I_{9}\,I_{3}^2 +
    213912\,I_{27}\,I_{6}\,I_{3}^3 - 128\,I_{27}\,I_{3}^5%
    )\,,
  \end{dmath*}
  \begin{dmath*}
    7\:\:
    \Phi_3(\betab_{22}) = %
    \text{-}\,{2^{135}\cdot 3^5} \:\:(%
    540\,I_{27}^2\,J_{12} - 4590\,I_{27}^2\,I_{12} - 151200\,I_{27}^2\,I_{6}^2
    + 4005\,I_{27}^2\,J_{9}\,I_{3} - 1683\,I_{27}^2\,I_{9}\,I_{3} -
    143010\,I_{27}^2\,I_{6}\,I_{3}^2 + 56\,I_{27}^2\,I_{3}^4%
    )\,.
  \end{dmath*}
\end{dgroup*}
Beside Klein's formula
\begin{small}
  $\Phi_3(\chib_{18}) = \text{-}\,2^{108} \:\:I_{27}^2$\,,
\end{small}
one finds a surprisingly compact expression for $\chib_{28}$\,,
\begin{displaymath}
  \Phi_3(\chib_{28}) = \text{-}\, 2^{171}\cdot 3^3 \:\:I_{27}^3 \, I_3\,.
\end{displaymath}

If we do not not pay attention, the rational coefficients of these formulas
tend to have prime factors greater than 7 in their denominators, especially
for the forms of higher weight.
We have eliminated all these ``bad primes'' using the relations that exist
between the Dixmier-Ohno invariants. It is also a good way to reduce
 the size of these expressions significantly. All in all, we gain a factor of 3
in the amount of memory to store the results (\textit{cf.} Table~\ref{tab:MFfromDO}).
\begin{table}[htbp]\centering
\begin{small}
  \begin{displaymath}
    \begin{array}[t]{c|l|r|r|}
      \text{\scriptsize Form} & \text{\scriptsize Content} & \text{\scriptsize Terms} & \text{\scriptsize Digits}\\\hline
      \alphab_{4}& \ 2^{20} \cdot 3^{3} \cdot 7& 6& 6 \\
      \alphab_{6}& \text{-}2^{28} \cdot 3^{4} \cdot 5^{\text{-}1} \cdot 7^{\text{-}1}& 19& 12 \\
      \alphab_{10}& \text{-}2^{44} \cdot 3^{3} \cdot 5^{\text{-}4} \cdot 7^{\text{-}2}& 98& 23 \\
      \alphab_{12}& \ 2^{75} \cdot 3& 4& 5 \\
      \alphab_{12p}& \ 2^{52} \cdot 3^{2} \cdot 5^{\text{-}4} \cdot 7^{\text{-}3}& 200& 26 \\
      \betab_{14}& \ 2^{81} \cdot 3^{4}& 11& 8 \\
      \alphab_{16}& \ 2^{66} \cdot 3^{4} \cdot 5^{\text{-}6} \cdot 7^{\text{-}3}& 703& 35 \\
      \betab_{16}& \ 2^{83} \cdot 3^{4} \cdot 5^{\text{-}2} \cdot 7^{\text{-}2}& 29& 16 \\
      \chib_{18}& \text{-}2^{108}& 1& 1 \\
      \alphab_{18}& \text{-}2^{80} \cdot 3^{3} \cdot 5^{\text{-}5} \cdot 7^{\text{-}2}& 813& 36 \\\hline
    \end{array}
    \ \ \ %
    \begin{array}[t]{c|l|r|r|}
      \text{\scriptsize Form} & \text{\scriptsize Content} & \text{\scriptsize Terms} & \text{\scriptsize Digits}\\\hline
      \alphab_{20}& \ 2^{75} \cdot 3^{2} \cdot 5^{\text{-}13} \cdot 7^{\text{-}6}& 1941& 52 \\
      \gammab_{20}& \ 2^{122} \cdot 3^{4}& 2& 3 \\
      \betab_{22}& \text{-}2^{135} \cdot 3^{5} \cdot 7^{\text{-}1}& 7& 6 \\
      \betab_{22p}& \ 2^{90} \cdot 3^{4} \cdot 5^{\text{-}14} \cdot 7^{\text{-}6}& 4000& 56 \\
      \alphab_{24}& \ 2^{89} \cdot 3^{2} \cdot 5^{\text{-}17} \cdot 7^{\text{-}7}& 6572& 67 \\
      \gammab_{24}& \ 2^{96} \cdot 3^{3} \cdot 5^{\text{-}14} \cdot 7^{\text{-}7}& 6585& 62 \\
      \gammab_{26}& \text{-}2^{105} \cdot 3^{5} \cdot 5^{\text{-}17} \cdot 7^{\text{-}7}& 10750& 67 \\
      \chib_{28}& \ 2^{171} \cdot 3^{3}& 1& 1 \\
      \alphab_{30}& \ 2^{109} \cdot 3 \cdot 5^{\text{-}21} \cdot 7^{\text{-}9}& 25630& 86\\\hline
    \end{array}
  \end{displaymath}
\end{small}\medskip
\caption{Polynomial expressions of the modular forms from Theorem~\ref{th:mfbasis}
  in terms of Dixmier-Ohno invariants: their content, their number of monomials,
  and the number of digits of the largest coefficient of their primitive part.}
\label{tab:MFfromDO}
\end{table}

\begin{remark}\label{rem:modforminDO}
  When we deal with the Jacobian of a curve with coefficients in $\Q$, what
  is a matter of few integer arithmetic operations to evaluate modular forms
  from invariants is a matter of high precision floating point arithmetic over
  the complex numbers with analytic computations of Riemann matrices. In practical
  calculations, such as the computations in Section~\ref{sec:relations}, it is
  thus much better to use the former, since a calculation that would take the
  order of the minute ultimately requires only a few milliseconds.
\end{remark}

\subsection{Invariants in terms of modular functions}
\label{sec:invar-terms-modul}

Conversely, we can look for expressions of a generating set of invariants in
terms of modular forms.  Using~\cite{Tsuyumine86,lorenzo19}, one obtains such
a result for invariants of binary octics. We focus here on the case of
Dixmier-Ohno invariants.

Since the locus of plane quartic over $\C$ such that $I_{27} \ne 0$
corresponds to the locus of non-hyperelliptic curve of genus $3$ and then to
principally polarized abelian threefolds $\C^3/(\tau \Z^3 + \Z^3)$ for which
$\chib_{18}(\tau) \ne 0$ \cite[Lem.~10,~11]{igusag3}, we see that any
invariant in $\bfI_3$ can be obtained as a quotient of a modular form by a
power of $I_{27}$.

\begin{proposition}\label{prop:phiinv}
  Let $I$ be a Dixmier-Ohno invariant of degree $3k$. There exist a polynomial
  $P_I$ in the modular forms from Theorem~\ref{th:mfbasis}, of weight
  $28\,k$, such that
  \begin{equation}\label{eq:8}
    I_{27}^{3k}\cdot I = \Phi_3(\,P_I(\,\alphab_4,\, \alphab_6, \ldots, \alphab_{30}\,)\,)\,.
  \end{equation}
\end{proposition}
\noindent
The first ones\footnote{We make available the list of these 13 polynomials
  at~\cite[file ``\texttt{SiegelMFtoDO.txt}'']{LR2019}.} are
\begin{dgroup*}[style={\footnotesize},spread={-2pt}]
  \begin{dmath*}
    2^{171}\cdot 3^3\:\:
    I_{27}^3 \: I_3 = \Phi_3(\text{-}\chib_{28})\,,
  \end{dmath*}
  \begin{dmath*}
    2^{144}\cdot 3^8\cdot 5\:\:
    I_{27}^6 \: I_6 = \Phi_3(\chib_{28}^2 - 2^4\cdot 3^2\: \chib_{18}^2\,\gammab_{20})\,,
  \end{dmath*}
  \begin{dmath*}
    2^{515} \cdot 3^{12} \cdot 5 \cdot 7^4\:\:
    I_{27}^9 \: I_9 =%
    \Phi_3(\:(-11735539200\,\alphab_{12}' - 2920548960\,\alphab_{12} -
    86929920\,\alphab_{6}^2 - 2027520\,\alphab_{4}^3)\:\chib_{18}^4+%
    (3259872\,\betab_{16}\,\betab_{14} - 4074840\,\gammab_{20}\,\alphab_{10} +
    21732480\,\gammab_{24}\,\alphab_{6} - 24837120\,\gammab_{26}\,\alphab_{4}
    + 137984\,\gammab_{20}\,\alphab_{6}\,\alphab_{4})\:\chib_{18}^3+%
    153856080\,\chib_{28}\,\gammab_{20}\:\chib_{18}^2%
    -1764735\,\chib_{28}^3\:) \,,
  \end{dmath*}
  \begin{dmath*}
    2^{515} \cdot 3^{12} \cdot 5^2 \cdot 7^4\:\:
    I_{27}^9 \: J_9 =%
    \Phi_3(\:(-30939148800\,\alphab_{12}' - 2200413600\,\alphab_{12} -
    229178880\,\alphab_{6}^2 - 5345280\,\alphab_{4}^3)\,\chib_{18}^4+%
    (8594208\,\betab_{16}\,\betab_{14} - 10742760\,\gammab_{20}\,\alphab_{10}
    + 57294720\,\gammab_{24}\,\alphab_{6} -
    65479680\,\gammab_{26}\,\alphab_{4} +
    363776\,\gammab_{20}\,\alphab_{6}\,\alphab_{4})\,\chib_{18}^3%
    +558376560\,\chib_{28}\,\gammab_{20}\,\chib_{18}^2%
    -5294205\,\chib_{28}^3\:) \,.
  \end{dmath*}
\end{dgroup*}
In this setting, one can also write $I_{27}^{27} \: I_ {27} =\Phi_3((2^{- 108} \,
\chib_ {18})^{14}) $.
\smallskip

Unlike the previous computations, one cannot obtain the above ones by a
direct application of the evaluation/interpolation strategy as the degrees
(and weights) are sometimes too large. For the invariant $I_{21}$, for
instance, one would potentially need to interpolate on a vector space of
modular forms of weight 196, which is huge (its dimension is $869\,945$).  The
trick is to proceed by steps and first look for expressions of a small power
of $I_ {27}$ by the desired invariant $I$, not only in terms of modular forms,
but also in terms of invariants $ I_{3k}$ of smaller degrees. For instance in
the case of $I_{21}$,
\begin{dgroup*}[style={\footnotesize},spread={-2pt}]
  \begin{dmath*}
    2^{63}\cdot3^{21}\cdot5^{21}\cdot7^{10}\cdot 11\:\: I_{27}\,I_{21} =
    2^{51}\cdot3^{15}\cdot5^{18}\cdot7^{9}\cdot 11 \:\:I_{27}\, ( -16156800\,J_{12}\,J_{9} +
    5680595070\,I_{12}\,J_{9} + 109296000\,J_{12}\,I_{9} - 3076972650\,I_{12}\,I_{9} -
    216169581600\,J_{15}\,I_{6} + 439538400\,I_{15}\,I_{6} - 770217033600\,J_{9}\,I_{6}^2 +
    2235454502400\,I_{9}\,I_{6}^2 + 8070768720\,J_{18}\,I_{3} - 622051920\,I_{18}\,I_{3} -
    3928070295\,J_{9}^2\,I_{3} + 1754339490\,J_{9}\,I_{9}\,I_{3} - 182964375\,I_{9}^2\,I_{3} +
    70135124400\,J_{12}\,I_{6}\,I_{3} - 611730004680\,I_{12}\,I_{6}\,I_{3} -
    18401013388800\,I_{6}^3\,I_{3} - 8799659820\,J_{15}\,I_{3}^2 + 1352865780\,I_{15}\,I_{3}^2 +
    237928085190\,J_{9}\,I_{6}\,I_{3}^2 - 56462733090\,I_{9}\,I_{6}\,I_{3}^2 + 294430290\,J_{12}\,I_{3}^3
    - 1980696900\,I_{12}\,I_{3}^3 - 4995876680760\,I_{6}^2\,I_{3}^3 - 65637369\,J_{9}\,I_{3}^4 -
    76264307\,I_{9}\,I_{3}^4 + 4016874680\,I_{6}\,I_{3}^5 )
    + 2^{12} \cdot 3^{6} \cdot 5^{3} \cdot 7\:\Phi_3(-19003712\,\betab_{16} - 10671360\,\alphab_{16} - 11116\,\alphab_{10}\,\alphab_{6 }-
    1844513\,\alphab_{12}\,\alphab_{4})\,.
  \end{dmath*}
\end{dgroup*}
Then, mechanically, through a sequence of substitutions of the invariants of
smaller degrees by their expression in terms of the modular forms, we arrive
to expressions for $ I_ {27}^{3k} \, I_ {3k} $ purely in terms of modular
forms. These formulas are very sparse, considering their weight (see
Table~\ref{tab:DOfromMF}).

\begin{table}[htbp]\centering
\begin{small}
  \begin{displaymath}
    \begin{array}[t]{c|l|r|r|}
      \text{\scriptsize DO inv.} & \text{\scriptsize Content} & \text{\scriptsize Terms} & \text{\scriptsize Digits}\\\hline
I_{3} & 2^{\text{-}171} \cdot 3^{\text{-}3} & 1 & 1\\
I_{6} & 2^{\text{-}344} \cdot 3^{\text{-}8} \cdot 5 & 2 & 3\\
I_{9} & 2^{\text{-}515} \cdot 3^{\text{-}12} \cdot 5 \cdot 7^{\text{-}4} & 11 & 11\\
J_{9} & 2^{\text{-}515} \cdot 3^{\text{-}12} \cdot 5^{\text{-}2} \cdot 7^{\text{-}4} & 11 & 11\\
I_{12} & 2^{\text{-}686} \cdot 3^{\text{-}16} \cdot 5^{\text{-}2} \cdot 7^{\text{-}4} & 13 & 13\\
J_{12} & 2^{\text{-}686} \cdot 3^{\text{-}16} \cdot 5 \cdot 7^{\text{-}3} & 14 & 13\\
I_{15} & 2^{\text{-}859} \cdot 3^{\text{-}20} \cdot 5^{\text{-}2} \cdot 7^{\text{-}4} \cdot 11^{\text{-}1}& 58 & 17\\
J_{15} & 2^{\text{-}859} \cdot 3^{\text{-}18} \cdot 5^{\text{-}3} \cdot 7^{\text{-}4} \cdot 11^{\text{-}1}& 58 & 17\\
I_{18} & 2^{\text{-}1030} \cdot 3^{\text{-}24} \cdot 5^{\text{-}2} \cdot 7^{\text{-}7} \cdot 11^{\text{-}2} \cdot 19^{\text{-}1}& 1321 & 237\\
J_{18} & 2^{\text{-}1030} \cdot 3^{\text{-}24} \cdot 5^{\text{-}3} \cdot 7^{\text{-}8} \cdot 11^{\text{-}3} \cdot 19^{\text{-}1}& 1321 & 238\\
I_{21} & 2^{\text{-}1202} \cdot 3^{\text{-}29} \cdot 5^{\text{-}6} \cdot 7^{\text{-}8} \cdot 11^{\text{-}3} \cdot 19^{\text{-}1}& 1382 & 242\\
J_{21} & 2^{\text{-}1201} \cdot 3^{\text{-}27} \cdot 5^{\text{-}6} \cdot 7^{\text{-}8} \cdot 11^{\text{-}3} \cdot 19^{\text{-}1}& 1382 & 242\\
I_{27} & 2^{\text{-}1512} & 1 & 1\\\hline
  \end{array}
  \end{displaymath}
\end{small}\medskip
\caption{Polynomial expressions of the Dixmier-Ohno in terms of the $19$
  generators from Theorem~\ref{th:mfbasis}: the content, the number of
  monomials, and the number of digits of the largest coefficient of the
  primitive parts.}
\label{tab:DOfromMF}
\end{table}

\begin{remark}
  It is not a coincidence that the power of $I_{27}$ is $3\,k$ in
  Equation~(\ref{eq:8}). Let us consider ternary quartics of the form
  $Q^2+p\,G$, where $p$ is a prime integer, $Q$ is a ternary quadratic form
  and $G$ is ternary quartic form. Generically, for all but $I_{27}$ the
  valuation of $p$ of the Dixmier-Ohno invariants of these forms is zero, and
  $\upsilon_p(I_{27})=14$\,. And, still generically, we have
  $\upsilon_p(\Phi_3(f))=3\,h/2$ where $f$ is any one of the Tsuyumine modular
  forms, and $h$ is its weight. Thus, if the equation
  $I_{27}^{\,\kappa}\cdot I = \Phi_3(\,P_I(\,\alphab_4,\, \alphab_6, \ldots,
  \alphab_{30})\,)$ is satisfied, the power $\kappa$ of $I_{27}$ must be such that
  the degrees agree, \textit{i.e.}  $27\,\kappa+3\,k = 3\,h$, and such
  that the valuations at $p$ are equal, \textit{i.e.} $14\,\kappa =
  3\,h/2$. This yields $\kappa = 3\,k$.
\end{remark}

\begin{remark}
  We are also able to eliminate the primes greater than $7$ in the
  denominators of the coefficients in these formulas using the relations that
  exist between  Siegel modular forms (\textit{cf.}
  Section~\ref{sec:relations}), with the notable exception of the primes $11$ and
  $19$ (\textit{cf.}  Table~\ref{tab:DOfromMF}).
  We suspect that the reason behind this difficulty is that, similarly to the
  prime $11$ (\textit{cf.} Remark~\ref{rem:mod11}), one cannot extend
  Theorem~\ref{th:mfbasis} \textit{mutatis mutandis} to characteristic
  $19$. Although we do not go further into the topic, it is possible to work
  directly in these characteristics and find specific formulas valid there.
\end{remark}

\newcommand{\etalchar}[1]{$^{#1}$}

\end{document}